\documentclass[10pt,journal,compsoc]{IEEEtran}
\pdfoutput=1

\ifCLASSOPTIONcompsoc
  \usepackage[nocompress]{cite}
\else
  \usepackage{cite}
\fi
\usepackage[pdftex]{graphicx}
\usepackage{amsmath}
\usepackage{algorithmic}
\usepackage{algorithm}

\usepackage{array}
\begin{document}
%
\title{A Numerical-based Parametric Error Analysis Method for Goldschmidt Floating Point Division}
%
%
%
%

\author{Binzhe~Yuan,
         Liangtao~Dai,
        and~Xin~Lou,~\IEEEmembership{Senior~Member,~IEEE}
\IEEEcompsocitemizethanks{\IEEEcompsocthanksitem The authors are with the School of Information Science and Technology, ShanghaiTech University, Shanghai 201210, China.\protect\\
}
}

\IEEEtitleabstractindextext{%
\begin{abstract}
This paper proposes a parametric error analysis method for Goldschmidt floating point division, which reveals how the errors of the intermediate results accumulate and propagate during the Goldschmidt iterations.
The analysis is developed by separating the error terms with and without convergence to zero, which are the key parts of the iterative approximate value. The proposed method leads to a state-of-the-art wordlength reduction for intermediate results during the Goldschmidt iterative computation. 
It enables at least half of the calculation precision reduction for the iterative factor to implement the rectangular multiplier in the divider through flexible numerical method, 
which can also be applied in the analysis of Goldschmidt iteration with iterative factors assigned under other ways for faster convergence.
Based on the proposed method, two proof-of-concept divider models with different configurations are developed, which are verified by more than 100 billion random test vectors to show the correctness and tightness of the proposed error analysis.
\end{abstract}

}

\maketitle

\IEEEdisplaynontitleabstractindextext

%
\IEEEpeerreviewmaketitle

\IEEEraisesectionheading{\section{Introduction}\label{sec:introduction}}

%
%
%
%
\IEEEPARstart{}{Floating} point divider is one of the fundamental building blocks in various processors.
For efficient designs, digit-recurrence and functional iteration are two main types of algorithms implemented in floating-point dividers of modern processors.
Iterative methods such as Newton-Raphson\cite{beame1986log} and Goldschmidt\cite{goldschmidt1964applications} algorithms usually have quadratic convergence, making them significantly faster than the digit-recurrence methods such as non-restoring and Sweeney-Robertson-Tocher (SRT)\cite{parhami2010computer} algorithms, which have linear convergence. Moreover, the Goldschmidt method turns the two multiplications in Newton's approach from dependent to parallelizable so as to halve the theoretical calculation latency. Therefore, Goldschmidt algorithms are more suitable for high speed floating-point divider design.

However, compared with digit-recurrence methods, functional iteration approaches such as Goldschmidt can only derive approximate quotient without exact information of the remainder. Current rounding methods\cite{cornea1999correctness, obermann1997division, schulte2007floating, kong2010rounding, piso2010variable} usually require an approximate quotient whose error is limited within a certain range, and a metric of its error with respect to the quotient of infinite precision is required, e.g. the remainder, to realize the correct rounding.

In practical circuit implementation of Goldschmidt algorithm, calculation errors are usually introduced by the truncation of wordlength of intermediate results. The complexity the iteration formula increases exponentially after the introduction of the aforementioned errors, making it hard to analyze and derive the accurate range of the final error to guarantee the correct rounding. Since all the errors come from finite wordlength, the most straightforward way to reduce the final error is to increase the intermediate calculation precision. But the hardware cost will increase correspondingly. As a result, it is crucial to develop an error analysis theory that can not only ensure correct rounding but also minimize circuit complexity. 

Existing error analysis methods usually have complicated form or certain limitations from the hardware implementation perspective. The theory that Goldschmidt proposed in \cite{goldschmidt1964applications} focuses mainly on the specific divider design in that work. Thereafter, another more general analysis method was proposed in \cite{krishnamurthy1970optimal}. But it allows only one multiplicand to be imprecise, limiting its applications. In the 1990s, dividers based on the Goldschmidt algorithm were applied in several commercial products of IBM\cite{schwarz1997cmos} and AMD\cite{oberman1999floating}. In one of these designs, i.e., the floating-point divider in AMD-K7, the designers adopt an error analysis method that over-estimates the intermediate error to ensure the correctness of the final results. Moreover, it mainly depends on a mechanically-checked formal proof\cite{russinoff1998mechanically}. Afterwards, a new parametric error analysis method was proposed in \cite{even2005parametric}, which formulates intermediate errors in a relative form and derives more elaborate and precise results. Moreover, the error analysis in \cite{kong2010rounding} presents a more concise way to approach the final error bound and the multiplications only needs $3$ extra bits of precision.

For Goldschmidt iteration itself, the work in \cite{schulte2007floating} indicates that the precision of intermediate calculation can be reduced under a certain rule, and rectangular multipliers whose two operands have different wordlengths can be adopted to improve the performance. It has also been shown that the speed of convergence can be faster than quadratic if the value of the iterative factor is assigned under a new way\cite{kong2009goldschmidt}.

In this paper, we propose a parametric error analysis method for Goldschmidt floating point division to reveal how each of the errors affects the final approximate result when the intermediate calculation is inaccurate. The analysis in this work is developed by separating the error terms with and without convergence to zero, which are the key parts of the iterative approximate value. In the derivative process, how the wordlength of the operands can be reduced will be shown naturally, which provides the theoretical foundation of taking advantage of rectangular multiplier in the divider's circuit implementation. With the proposed theory, principles for determining each of the wordlengths can be discussed under the assumptions of initial reciprocal approximation and final rounding method.

This paper is organized as follows. Section \ref{sec:background} introduces the fundamental of the original Goldschmidt algorithm and its implements.
In Section \ref{sec:error_analysis}, the modified Goldschmidt algorithm with intermediate calculation errors is proposed, followed by the error analysis and two important error terms. The final error bound is derived in Section \ref{sec:error_bound}, and two practical examples of wordlength selection are presented in Section \ref{sec:implementation}.
Finally, several criticality problems of the proposed method and the verification results for the two divider examples are discussed in Section \ref{sec:discussion}.

\section{Background}\label{sec:background}

The Goldschmidt algorithm represents the division of $Q=A/B$ in a fraction form, and multiplies both the numerator and denominator with one of a series of specific iterative factors for correction in each iteration. When the denominator approaches to $1$, the numerator will correspondingly approach to the quotient $Q$. Let us denote the numerator, denominator and iterative factor of each iteration as $N_i$, $D_i$ and $F_i$. The procedure of the algorithm is described in Algorithm \ref{original_GS}, where $Q_A$ indicates the approximate quotient value generated by the Goldschmidt iteration. For the value of $F_i$, the Goldschmidt algorithm proposes a method which is easy to implement by circuit, e.g. $F_i=2-D_i$. The initial iterative factor is assigned with the value of $rec$, which is the approximate reciprocal of the divisor $B$ with relative error less than $1$. 

\begin{algorithm}[t]
\caption{Goldschmidt Division of $Q=A/B$}
\label{original_GS}
\begin{algorithmic}
\STATE \textbf{Initialization:} $N_{-1}=A$, $D_{-1}=B$ and $F_{-1}=rec$
\FOR{$i=0$ to $k$}
    \STATE $N_{i}=N_{i-1}\times F_{i-1}$
    \STATE $D_{i}=D_{i-1}\times F_{i-1}$
    \STATE $F_{i}=2-D_{i-1}$
\ENDFOR
\STATE $Q_{A}=N_{k}$
\end{algorithmic}
\end{algorithm}

As a basis for the subsequent error analysis, it is necessary to first discuss the algorithms for the reciprocal look-up table and rounding method adopted in this paper. The bipartite table proposed in \cite{sarma1995faithful, schulte1997symmetric} only needs two small memory blocks and a low bit-width subtractor to implement the reciprocal table with relatively high precision, making it faster than the table based on simple interpolation\cite{farmwald1981design, nakano1987method}. 
In this work, the bipartite table is implemented in a more flexible way, i.e., parameter tuning of the input and output data-width of the bipartite table is of higher freedom.
The error analysis of a particular table is solved by traversing all the elements in the table, calculating each local maximum relative error and deriving the global value. Due to the flexible parameter tuning, a series of table configurations with corresponding maximum relative error that has quite continuous values can be derived. These particular maximum relative errors help to find the table configuration as ideally as possible. 

As for the rounding method, the algorithm mentioned in \cite{kong2010rounding} claims that it doubles the error tolerance and implements a correct rounding under the IEEE 754 standard, $\pm 2^{-(n+1)}$ is the bound for the absolute error of the approximate quotient. Here, $n$ refers to the number of bits for the mantissa in the target floating-point number where the quotient of infinite precision has a format of $1.b_1b_2b_3\cdots$. 

\section{Error Analysis}\label{sec:error_analysis}
For circuit implementation of Goldschmidt division, the three iterative variables are usually not precisely computed for complexity consideration, i.e., the original products from the multiplier need to be rounded/truncated before the next iteration. Correspondingly, there are three different types of errors introduced by rounding/truncation: (i) the absolute error on $N_i$, denoted by $n_i$; (ii) the absolute error on $D_i$, denoted by $d_i$ and (iii) the absolute error on $F_i$, denoted by $f_i$. Moreover, for the convenience of hardware implementation, truncation is used in this work, meaning that the errors $n_i,d_i,f_i$ are all non-negative. In this section, we first analyze the errors of three iterative variables independently, and propose a complete error analysis by considering all the three error terms.

\subsection{Numerator Error}
Consider the case that the truncation error appears only in the calculation of $N_i$, and if $N_i$ remains $k$ fractional bits after truncation, then we have $0\leq n_i<2^{-k}=ulp$. For iteration 0, the Goldschmidt equation with truncation error is given by
\begin{IEEEeqnarray}{rCl}
\frac{N_0}{D_0}&=&\frac{A \cdot rec -n_0}{B \cdot rec}\nonumber\\
               &=&\frac{Q(1-\epsilon_0)-n_0}{1-\epsilon_0}\label{eqn:N_it0}\\
            F_0&=&2-D_0=1+\epsilon_0\nonumber
\end{IEEEeqnarray}
Note that the term $Q(1-\epsilon_0)$ in the numerator of (\ref{eqn:N_it0}) will keep converging to $Q$ and the denominator $1-\epsilon_0$ will converge to $1$, which are the same as that in the accurate calculation. Another term $n_0$ will not converge since it dose not contain the same converging sub-term $(1-\epsilon_0)$.
The iteration equations for iteration 1 and iteration 2 can be derived as
\begin{IEEEeqnarray*}{rCl}
\frac{N_1}{D_1}&=&\frac{N_0 \cdot F_0 -n_1}{D_0 \cdot F_0}\nonumber\\
               &=&\frac{Q(1-\epsilon_0^2)-n_0F_0-n_1}{1-\epsilon_0^2}\\
            F_1&=&2-D_1=1+\epsilon_0^2
\end{IEEEeqnarray*}
and
\begin{IEEEeqnarray*}{rCl}
\frac{N_2}{D_2}&=&\frac{N_1 \cdot F_1 -n_2}{D_1 \cdot F_1}\nonumber\\
              &=&\frac{Q(1-\epsilon_0^4)-n_0F_0F_1-n_1F_1-n_2}{1-\epsilon_0^4}\\
            F_2&=&2-D_2=1+\epsilon_0^4
\end{IEEEeqnarray*}
respectively. If the initial error $\epsilon_0$ is small enough, then the iterative factor $F_i$ is very close to $1$, meaning that being multiplied by the iterative factor $F_i$ in each iteration has little effect on the terms that contain $n_i$ in the numerator. According to this, a significant approximation and inference can be made:
\newtheorem{lemma}{Lemma}[section]
\begin{lemma}
Given a very small initial error $\epsilon_0$, the effect of the multiplication of $n_i$ terms and the iterative factors $F_i$, which is very close to $1$, can be ignored conditionally, i.e. $n_{i}(1+\epsilon_0^{2^i})(1+\epsilon_0^{2^{i+1}})\cdots \approx n_i$. This will be discussed in detail later.
\label{lem:ignore_Fi_on_N}
\end{lemma}
\newtheorem{corollary}{Corollary}[section]
\begin{corollary}
By representing the upper bound of all $n_i$ with a unified $n$ and assuming that the number of the $n_i$ terms in the numerator is $m$, a simple upper bound for the sum of all $n_i$ terms in the numerator can be represented as
\begin{equation*}
    \sum (terms\;with\;n_i) <(m+1)\cdot n.
\end{equation*}
\label{cor:ignore_Fi_on_N}
\end{corollary}
As a result, the approximate expression of $\frac{N_2}{D_2}$ can be expressed as
\begin{equation*}
    \frac{N_2}{D_2} \approx \frac{Q(1-\epsilon_0^4)-n_0-n_1-n _2}{1-\epsilon_0^4}.
\end{equation*}
This expression reveals that the error caused by the inaccurate calculation of the numerator will not converge to $0$, i.e. the $n_i$ error is not self-correcting. If $N_2$ is the iterative result, the lower bound of its error can be simply constrained by
\begin{IEEEeqnarray*}{rCl}
    N_2-Q&=&-Q\epsilon_0^4-n_0F_0F_1-n_1F_1-n_2\nonumber\\
         &>&-Q\epsilon_0^4-(3+1)\cdot n\nonumber\\
         &>&-Q\epsilon_0^4-(3+1)\cdot ulp.
\end{IEEEeqnarray*}

\subsection{Denominator Error}
Consider the case that the truncation error appears only in the calculation of $D_i$. Similarly, let $D_i$ remains $k$ fractional bits after truncation, then the equation for iteration 0 becomes
\begin{IEEEeqnarray*}{rCl}
    \frac{N_0}{D_0}&=&\frac{A \cdot rec}{B \cdot rec-d_0}\nonumber\\
                   &=&\frac{Q(1-\epsilon_0)}{1-\epsilon_0-d_0},\\
                F_0&=&2-D_0=1+\epsilon_0+d_0.
\end{IEEEeqnarray*}
In order to separate the parts that with and without convergence, let us define $\epsilon_0^{\prime}=\epsilon_0+d_0$ to track the term in the numerator that converges along with the denominator, the original equation can then be rewritten as
\begin{IEEEeqnarray*}{rCl}
\frac{N_0}{D_0}&=&\frac{Q(1-\epsilon_0^{\prime})+Qd_0}{1-\epsilon_0^{\prime}},\\
            F_0&=&1+\epsilon_0^{\prime},
\end{IEEEeqnarray*}
whose form is similar to (\ref{eqn:N_it0}). Therefore, the $Q(1-\epsilon_0^{\prime})$ term will keep the property of converging with the denominator while the $Qd_0$ term will not.

Likewise, for iteration 1 and 2, let us define $\epsilon_i^{\prime}=\epsilon_{i-1}^{\prime 2}+d_i$ and $F_i=1+\epsilon_i^{\prime}$, the results of the iterations are as follows
\begin{IEEEeqnarray*}{rCl}
    \frac{N_1}{D_1}&=&\frac{Q(1-\epsilon_1^{\prime})+Qd_0F_0+Qd_1}{1-\epsilon_1^{\prime}},\\
    \frac{N_2}{D_2}&=&\frac{Q(1-\epsilon_2^{\prime})+Qd_0F_0F_1+Qd_1F_1+Qd_2}{1-\epsilon_2^{\prime}}.
\end{IEEEeqnarray*}
Therefore, similar approximation and inference as in the numerator error analysis can be made:
\begin{lemma}
The effect of the multiplication of $Qd_i$ terms and the iterative factors $F_i$, which is very close to $1$, can be ignored conditionally, i.e. $Qd_i(1+\epsilon_{i}^{\prime})(1+\epsilon_{i+1}^{\prime})\cdots \approx Qd_i.$
\label{lem:ignore_Fi_on_D}
\end{lemma}
\begin{corollary}
By representing the upper bound of all $d_i$ with a unified $d$ and assuming that the number of the $Qd_i$ terms in the numerator is $m$, then a simplified upper bound for the sum of all $Qd_i$ terms in the numerator can be represented as
\begin{equation*}
    \sum (terms\;with\;Qd_i) <(Qm+1)\cdot d.
\end{equation*}
\label{cor:ignore_Fi_on_D}
\end{corollary}

Similar to the property of $n_i$ error, the $d_i$ error is also not self-correcting. The difference between these two errors is that the $d_i$ error turns the convergent term $(1-\epsilon_i)$ to $(1-\epsilon_i^{\prime})$, which may affect the speed of convergence for the subsequent iterations. However, except in the last iteration where the upper bound of $\epsilon_i$ is with similar or less magnitude comparing to the $d_i$ error itself, the difference between $\epsilon_i$ and $\epsilon_i^{\prime}$ is small enough compared to their upper bound in every case, such that it does not influence the trend of convergence greatly. If $N_2$ is the iterative result, the upper bound of its error can be constrained by
\begin{IEEEeqnarray*}{rCl}
    N_2-Q&=&-Q\epsilon_1^{\prime 2}+Qd_0F_0F_1+Qd_1F_1\nonumber\\
         &<&-Q\epsilon_1^{\prime 2}+(2Q+1)\cdot d\nonumber\\
         &<&(2Q+1)\cdot ulp.
\end{IEEEeqnarray*}

\subsection{Iterative Factor Error}
The iterative factor $F_i$ is computed by the subtracting $D_i$ from a constant $2$. Commonly, for circuit implementation, subtraction is implemented using two's complement, which requires to take the inverse of all bits and add $1$. This addition may increase the circuit and timing complexity. Therefore, if the error introduced is controllable, the iterative factor can be computed using one's complement instead, where the $f_i$ error introduced to $F_i$ is exactly $2^{-k}$, equaling to $ulp$.

If only the $f_i$ error appears in the iteration, the equation for iteration 0 becomes
\begin{IEEEeqnarray*}{rCl}
    \frac{N_0}{D_0}&=&\frac{A \cdot rec}{B \cdot rec}
                   =\frac{Q(1-\epsilon_0)}{1-\epsilon_0},\\
                F_0&=&2-D_0-f_0=1+\epsilon_0-f_0.
\end{IEEEeqnarray*}
Similarly, for iteration 1, we have
\begin{equation*}
    \frac{N_1}{D_1}=\frac{N_0\cdot F_1}{D_0 \cdot F_1}
                   =\frac{Q(1-\epsilon_0)(1+\epsilon_0-f_0)}{(1-\epsilon_0)(1+\epsilon_0-f_0)}.
\end{equation*}
Let us define $1-\epsilon_1=(1-\epsilon_0)(1+\epsilon_0-f_0)$, the result of iteration 1 can be rewritten as
\begin{IEEEeqnarray*}{rCl}
    \frac{N_1}{D_1}&=&\frac{Q(1-\epsilon_1)}{(1-\epsilon_1)},\\
                F_1&=&2-D_1-f_1=1+\epsilon_1-f_1,
\end{IEEEeqnarray*}
which shows that the $f_i$ error will not deprive the property that the numerator converges along with the denominator, i.e., the $f_i$ error is self-correcting. Consequently, the error caused by $f_i$ can be corrected by the subsequent iterations unless it is the last one. In the last iteration, $f_i$ will leave error that is not correctable in the numerator just like the $d_i$ error.

Since $f_i$ only leads to correctable error in iterations except the last one, the $f_i$ error which is far greater than $ulp$ can also be corrected. Consider a numerical example, let us randomly set $|\epsilon_0|=2^{-17.3}$. For accurate calculation, the error of the numerator will be $2^{-17.3\times 2}=2^{-34.6}$. However, if the iterative factor $F_i$ is truncated after one's complement and remains $35$ fractional bits, i.e. $f_0\leq 2^{-35}$, the numerator after iteration 1 is given by
\begin{equation*}
    N_1=Q(1-\epsilon_0)(1+\epsilon_0-f_0)=Q(1-2^{-17.3\times 2}-f_0(1-\epsilon_0)).
\end{equation*}
Therefore, the maximum error of the numerator will be less than
\begin{equation*}
    2^{-17.3\times 2}+2^{-35}\times (1+|\epsilon_0|)\approx 2^{-33.78},
\end{equation*}
which indicates a precision loss of not greater than $1$ bit comparing to $2^{-34.6}$. As a result, iterative factor of full precision is not necessary in the early stages of iteration. This explains the idea of reducing the wordlength of $F_i$ in the certain way mentioned by \cite{schulte2007floating}. Meanwhile, the iterative factor approaches to $1$ as the iteration carries on, the continuous digits of $0$ or $1$ after the binary point of $F_i$ can be ignored. Data-width reduction of multipliers in this work takes advantage of both the above properties.

\subsection{Complete Error Analysis}
In this subsection, we discuss the complete error analysis by considering all the three aforementioned different types of errors. The following equation describes the result of iteration 0:
\begin{IEEEeqnarray*}{rCl}
    \frac{N_0}{D_0}&=&\frac{A\cdot rec-n_0}{B\cdot rec-d_0}\\
                   &=&\frac{Q(1-\epsilon_0)-n_0}{1-\epsilon_0-d_0}\\
                   &=&\frac{Q(1-\epsilon_0^{\prime})+Qd_0-n_0}{1-\epsilon_0^{\prime}},\\
                F_0&=&2-D_0-f_0=1+\epsilon_0^{\prime}-f_0
\end{IEEEeqnarray*}
For iteration 1, let us define $D_0\cdot F_0=(1-\epsilon_0^{\prime})(1+\epsilon_0^{\prime}-f_0)=1-\epsilon_1$. The accurate value of $\epsilon_1$ should be calculated in numerical way to measure the degree of convergence. According to the property of Lemma \ref{lem:ignore_Fi_on_N} and \ref{lem:ignore_Fi_on_D}, the result of iteration 1 can be represented approximately by
\begin{IEEEeqnarray*}{rCl}
    \frac{N_1}{D_1}&=&\frac{N_0\cdot F_0-n_1}{D_0\cdot F_0-d_1}\nonumber\\
                   &\approx& \frac{Q(1-\epsilon_1)+Qd_0-n_0-n_1}{1-\epsilon_1-d_1}\nonumber\\
                   &=&\frac{Q(1-\epsilon_1^{\prime})+Q(d_0+d_1)-(n_0+n_1)}{1-\epsilon_1^{\prime}},\\
                F_1&=&2-D_1-f_1=1+\epsilon_1^{\prime}-f_1.
\end{IEEEeqnarray*}
Similarly, the approximate result of iteration 2 is given by
\begin{IEEEeqnarray*}{rCl}
    \frac{N_2}{D_2}&=&\frac{N_1\cdot F_1-n_2}{D_1\cdot F_1-d_2}\nonumber\\
                   &\approx& \frac{Q(1-\epsilon_2)+Qd_0-n_0+Qd_1-n_1-n_2}{1-\epsilon_2-d_2}\nonumber\\
                   &=&\frac{Q(1-\epsilon_2^{\prime})+Q(d_0+d_1+d_2)-(n_0+n_1+n_2)}{1-\epsilon_2^{\prime}},\\
                F_2&=&2-D_2-f_2=1+\epsilon_2^{\prime}-f_2.
\end{IEEEeqnarray*}
If iteration 3 is the last iteration, the approximate expression of $N_3$ can be derived in the same manner as
\begin{IEEEeqnarray}{rCl}
    N_3&=&N_2\cdot F_2-n_3\nonumber\\
       &\approx& Q(1-\epsilon_2^{\prime 2})-Qf_2(1-\epsilon_2^{\prime})+\sum_{i=0}^{2}Qd_i-\sum_{i=0}^{3}n_i.
       \label{eqn:complete_N3}
\end{IEEEeqnarray}
Generally, the $f_i$ error of the last iteration cannot be corrected, such that this $f_i$, denoted by $f_{last}$, should be as small as possible. If the last iteration factor is calculated using one's complement, then $f_{last}$ equals to $2^{-k}$ exactly, i.e. $ulp$.

Equation (\ref{eqn:complete_N3}) of $N_3$ shows how the final error is determined by the $n_i$, $d_i$ and $f_i$ errors approximately, and a strict error bound of it can be derived by Corollary \ref{cor:ignore_Fi_on_N} and \ref{cor:ignore_Fi_on_D}. For the upper bound, we replace all $n_i$ by $0$ and all $d_i$ by their corresponding upper bound in the equation to maximize the positive error as
\begin{IEEEeqnarray*}{rCl}
N_{3max}-Q&\approx& -Q\epsilon_2^{\prime 2}-ulp\cdot Q(1-\epsilon_2^{\prime})+3Qd\nonumber\\
          &\leq& -ulp\cdot Q(1-\epsilon_2^{\prime})+3Qd\nonumber\\
          &\approx& -ulp\cdot Q+3Q\cdot ulp\nonumber\\
          &=&2Q\cdot ulp\nonumber\\
          &<&2Q\cdot ulp+ ulp.
\end{IEEEeqnarray*}
To derive the lower bound, we set all $d_i$ to $0$ and all $n_i$ to their corresponding upper bound to minimize the negative error as
\begin{IEEEeqnarray*}{rCl}
N_{3min}-Q &\approx& -Q\epsilon_2^{\prime 2}-ulp\cdot Q(1-\epsilon_2^{\prime})-4n\nonumber\\
           &\approx& -Q\epsilon_2^{\prime 2} -Q\cdot ulp-4ulp\nonumber\\
           &>&       -Q\epsilon_2^{\prime 2} -Q\cdot ulp-4ulp -ulp.
\end{IEEEeqnarray*}
If the initial precision of the iteration is high enough to ensure the magnitude of the $-Q\epsilon_2^{\prime 2}$ term to be small enough, then the final error bound can be simplified as
\begin{equation*}
    (-Q\cdot ulp-4ulp-ulp \; ,\;2Q\cdot ulp+ulp).
\end{equation*}

\section{Error Bound}\label{sec:error_bound}
In the previous section, the properties of the three different types of intermediate errors have been discussed. The numerator error $n_i$ and denominator error $d_i$ lead to non-correctable error in the final quotient approximation. Besides, the error caused by the last iterative factor error $f_{last}$ will accumulate in the final result. Therefore, the terms that contain $n_i$, $Qd_i$ and $f_{last}$ in the final result can be classified into accumulative error term (AET). On the other hand, for the term that resembles $-Q\epsilon_{k-1}^{\prime 2}$, since it converges to $0$ as the iteration carries on, it can be classified into convergent error term (CET).

\subsection{Accumulative Error Term}\label{sec:accumulative_error}
For a Goldschmidt algorithm with $k$-times of iterations, its approximate accumulative error term (AAET) is given by
\begin{equation}
    \text{AAET}(k)=Q\cdot\sum_{i=0}^{k-1}{d_i}-\sum_{i=0}^{k}{n_i}-Q\cdot f_{last},
\label{eqn:accumulative_error}
\end{equation}
where $n_i,d_i\in [0,ulp)$, $f_{last}\in[ulp,2^{-b_{last}}]$ and $Q\in(\frac{1}{2},2)$.
It is worth noting that $n_i$ and $d_i$ refer to the truncation errors here. But as long as the corresponding error range is specified, this expression can be modified to adapt to the cases of other rounding modes. For $f_{last}$ error, the wordlength reduction operation should not be applied to the last iterative factor $F_{last}$ in the view of minimizing the final error bound. As we have discussed, one's complement could be an appropriate way to calculate the iterative factor, and under this assumption, $f_{last}$ equals to $ulp$. However, if there is extra data-width truncation after one's complement, then the upper bound of $f_{last}$ will correspond to another new value $2^{-b_{last}}$, where the $b_{last}$ here refers to the remaining fractional bits after truncation. While for other rounding modes applied to the iterative factor, the aforementioned expression is still adaptive if the error range of $f_{last}$ is modified correspondingly. For $Q$, the range of the two operands in this paper is $[1,2)$, meaning that the range of $Q$ is $(\frac{1}{2},2)$. This expression keeps adaptive to operands with other ranges as long as the corresponding range of $Q$ is specified.

Expression (\ref{eqn:accumulative_error}) is an approximate representation of its analytical expression, where approximations used here are based on the properties in Lemma \ref{lem:ignore_Fi_on_N} and \ref{lem:ignore_Fi_on_D}. To understand the effects of these approximations more intuitively, let us consider a counterexample with the precision of the initial look-up table less than $2^{-6}$. Under this assumption, after as many as two times of iterations, the error of the numerator can only converge to at most $2^{-24}$ ideally, which is even critical for correct rounding for single precision. For this reason, the initial iterative precision should be high enough in most cases for fast convergence. Without loss of generality, we assume that the initial iterative precision is beyond $2^{-6}$ in the following discussion.

With the assumption of $\epsilon_0\geq2^{-6}$, a certain upper bound for the factors which are ignored by Lemma \ref{lem:ignore_Fi_on_N}, and \ref{lem:ignore_Fi_on_D} in front of the terms that contain $n_i$ and $Qd_i$ in the numerator can be derived by calculating a certain strict upper bound for the iterative factors. The following numerical analysis result for Goldschmidt algorithm with 3 times of iterations ($f_{last}$ is not included since it does not multiply with any $F_i$) will show the difference between AAET and the true value:

\begin{IEEEeqnarray*}{rCl}
    \text{AET}(3)&=&F_0F_1F_2(Qd_0-n_0)\\
                   &&+\>F_1F_2(Qd_1-n_1)\\
                   &&+\>F_2(Qd_2-n_2)\\
                   &&-\>n_3,
\end{IEEEeqnarray*}
where
\begin{equation*}
    \lceil F_2 \rceil < \lceil F_1 \rceil < \lceil F_0 \rceil
\end{equation*}
and
\begin{equation*}
    \lceil F_0 \rceil=1+\lceil |\epsilon_0| \rceil+\lceil d_0 \rceil-\lfloor f_0 \rfloor\approx1+2^{-6}.
\end{equation*}
Then we have
\begin{IEEEeqnarray*}{rCl}
    F_0F_1F_2&<&F_0^3\approx 1.0476,\\
       F_1F_2&<&F_0^2\approx 1.0315,\\
          F_2&<&F_0\approx 1.0156,
\end{IEEEeqnarray*}
which can be used to derive a numerical rigorous bound for AET.

\begin{table}[t]
   \renewcommand{\arraystretch}{1.2}
  \setlength{\tabcolsep}{17pt}
    \centering
    \caption{The error of estimated AET through AAET.}
    \label{tab:AAETvsAET}
    \begin{tabular}{|c|c|c|}
    \hline
       Bound         &  AAET         & Numerical result\\
    \hline
    Upper bound & $+3Q\cdot ulp$ & $+3.0947Q\cdot ulp$\\
    \hline
    Lower bound & $-4\cdot ulp$  & $-4.0947\cdot ulp$\\
    \hline
    \end{tabular}
\end{table}

According to the data in Table \ref{tab:AAETvsAET}, in spite of scaling the factors $F_0F_1F_2$, $F_1F_2$ and $F_2$ for simplicity of calculation, the derived numerical range which is able to constrain the error of the final result strictly is close to the approximate range. This difference will be even smaller with a higher initial iterative precision. Consequently, the two lemmas and the corresponding corollaries are reasonable, when the error bound derived by AAET is far from the error tolerance of the actual design. It is appropriate to get a rough error bound by simply adding and subtracting one $ulp$ on AAET without careful numerical calculation just as Corollary \ref{cor:ignore_Fi_on_N} and \ref{cor:ignore_Fi_on_D}. Nevertheless, if the error bound derived by AAET is close to the error tolerance, e.g. the redundancy for AET is $3.1Q\cdot ulp$ which approaches the rough bound derived directly from AAET of $3Q\cdot ulp$. In this case, although the rough error bound is already beyond the error tolerance, a numerical error bound which is not only within this error tolerance but also able to constrain the error of the final result can be found by numerical method. In conclusion, the aforementioned analysis does not overestimate the accumulated error and can measure the magnitude of the error accurate enough to derive a tight error bound.

\subsection{Convergent Error Term}\label{sec:convergence_error}
We have discussed that the exact value and the trend of convergence for this error term depend on both the initial iterative precision $\epsilon_0$ and the $d_i$, $f_i$ errors of each iteration. Similarly, for $k$-times of iterations, the numerical value of the CET is given by
\begin{equation}
    \label{eqn:convergence_error}
    \text{CET}(i)=-Q\epsilon_{i-1}^{\prime 2},\quad 1\leq i\leq k,
\end{equation}
where
\begin{align*}
{\epsilon_i^{\prime}}&=
\begin{cases}
    \epsilon_i+d_i,& i=0,\nonumber\\
    \epsilon_{i-1}^{\prime 2}+(1-\epsilon_{i-1}^{\prime})f_{i-1}+d_i,& 1\leq i\leq k-1.
\end{cases}
\end{align*}
Despite the $\epsilon_{i-1}^{\prime}$ term is an important component of its corresponding iterative factor $F_i$ and CET $-Q\epsilon_{i-1}^{\prime 2}$, according to the discussion in Section \ref{sec:accumulative_error}, AAET derived by ignoring the effect of multiplying with $F_i$ is slightly different from the precise numerical result. Therefore, the value of $\epsilon_{i-1}^{\prime}$ that leads to a greater or less CET has little effect on AET, meaning that the two types of error terms can be treated and maximized independently to derive the complete error bound.

To maximize the magnitude of CET, the factor $\epsilon_{i}^{\prime}$ in each iteration should approach to its maximum. For iteration 0, we assign $\epsilon_{0}$ with the maximum relative error of the look-up table and set $d_0$ to $ulp$ to derive a strict and tight bound for $\epsilon_{0}^{\prime}$, i.e. $|\epsilon_0^{\prime}|\leq |\epsilon_0|+ulp$. Similar method can be applied to the following analysis, e.g. $|\epsilon_1^{\prime}|\leq \epsilon_0^{\prime 2}+(1+|\epsilon_0^{\prime}|)|f_0|+ulp$. It is worth noting that the absolute value which simplifies the calculation in the aforementioned inequation may lead to a slightly greater $\epsilon_{i}^{\prime}$, but the scaling is small enough to ensure a tight bound. Moreover, all calculations in (\ref{eqn:convergence_error}) are numerical, which makes it adaptive to all other methods to compute the iterative factor as long as the new iterative relationship and the range of $f_i$ are specified, e.g. the method in \cite{kong2009goldschmidt}.

Usually, there are three types of the degree of convergence for CET: $\lceil |\text{CET}|\rceil \gg \lceil |\text{AET}|\rceil$, $\lceil |\text{CET}|\rceil = \lceil |\text{AET}|\rceil$ or $\lceil |\text{CET}|\rceil \ll \lceil |\text{AET}|\rceil$ in order of magnitude. The first case happens in the early stage of iteration and CET dominates the error bound. The remaining two cases are in the last iteration. For the second case, tight error bound of CET can be calculated in numerical way and added to AET. While for the last case, how much CET effects on the final error is similar to that of ignoring the iterative factors in Lemma \ref{lem:ignore_Fi_on_N} and \ref{lem:ignore_Fi_on_D}. Therefore, the error bound derived by Corollary \ref{cor:ignore_Fi_on_N} and \ref{cor:ignore_Fi_on_D} can constrain the final error well.

Equation (\ref{eqn:convergence_error}) also implies that if $\epsilon_{i-1}^{\prime 2}$ and $f_{i-1}$ are close in order of magnitude, the $f_{i-1}$ error will significantly influences the value of $\epsilon_i^{\prime}$. In other words, smaller $f_{i-1}$ influences $\epsilon_i^{\prime}$ less, while the marginal benefit of decreasing $f_{i-1}$ by increasing the precision of $F_{i-1}$ diminishes. Additionally, the derived tight error bound of $\epsilon_i^{\prime}$ helps to determine the number of following bits after the binary point of $F_i$ that can be omitted in the multiplications, i.e. if $\lfloor log_{\frac{1}{2}}{|\epsilon_i^{\prime}|} \rfloor=e$, $F_i$ has a form of
\begin{equation*}
{F_i}=
\begin{cases}
    1.000\cdots00b_{e+1}b_{e+2}b_{e+3}\cdots,& 0<\epsilon_i^{\prime}< 2^{-e},\\
    0.111\cdots11b_{e+1}b_{e+2}b_{e+3}\cdots,& -2^{-e}< \epsilon_i^{\prime}<0.
\end{cases}
\end{equation*}
The omitting of the continuous $e-1$ bits of $1$ or $0$ when computing $N_i$($D_i$) can be realized by transforming the calculation to the multiply-add form as 
\begin{IEEEeqnarray}{rCl}
    \label{eqn:mult_to_FMA}
    N_i\times F_i&=&N_i+N_i\times(F_i-1)\nonumber\\
                 &=&N_i+(N_i\times F_i[e:b_{f_i}])\times2^{-e}.
\end{IEEEeqnarray}
Here, $b_{f_i}$ in (\ref{eqn:mult_to_FMA}) refers to the number of the remaining fractional bits after truncation to $F_i$, and $[ulp,2^{-b_{f_i}}]$ corresponds to the range of the respective $f_i$. Moreover, the multiplication in (\ref{eqn:mult_to_FMA}) is signed for $\epsilon_i^{\prime}$ with possibly positive or negative value. This case happens only in iteration 1 under the assumptions in this paper and for the following iterations, $\epsilon_i^{\prime}$ stays positive rigorously such that the sign bit, i.e.,  the $k$-th bit, can be further removed potentially. The above-mentioned two ways to reduce the wordlength of $F_i$ provide rigorous theoretical direction to the usage of rectangular multiplier (long operand of $N_i$ or $D_i$, short operand of $F_i$) in the divider\cite{schulte2007floating}.

\section{Implementation}\label{sec:implementation}
In Section \ref{sec:error_bound}, strategies of the derivation of reasonable numerical error bound have been discussed. In this section, several examples of practical divider designs are presented to validate the proposed error analysis method.
Recall that there are four different basic rounding modes in the IEEE 754 standard. The rounding algorithm in \cite{kong2010rounding} only requires a constrain for the error range of the approximate quotient. Therefore, for simplicity, the verification is realized by verifying whether the absolute error of the approximate quotient is within the range of $\pm 2^{-(n+1)}$ or not.

\subsection{Three-Stage Divider for SP, DP and EP}\label{sec:stage-3}
This example is a Goldschmidt divider design with 3 iterations, which is capable of generating approximate quotients with single, double and extend precision, i.e., SP, DP and EP, respectively. 
First of all, to simplify the rounding stage, a pretreatment of left-shift for $1$ bit is applied to the dividend $A$ if it is less than the divisor $B$. This operation ensures all the quotients of infinite precision to have the same format as required in \cite{kong2010rounding}, i.e. the range of $Q$ turns from $(\frac{1}{2},2)$ to $[1,2)$ such that the bounds of absolute error for the three types of precision are $2^{-24},\;2^{-53}$ and $2^{-64}$, respectively.

Secondly, which type of error term that dominates the absolute error in each iteration should be analyzed. Consider a rough analysis, for the three different precision types, look-up tables with maximum relative error of at most $2^{-12}$, $2^{-13.35}$ and $2^{-8}$ are required respectively. As a result, CET is the bottleneck of precision for the DP case. Meanwhile, for the EP case, if the requirement of absolute error for the DP case can be met after iteration 2, then the weight of CET in the absolute error for the EP case must be extremely small because of its quasi-quadratic convergence. Therefore, AET is the bottleneck of precision for the EP case while it is not for the DP case.

Next, assume that all the $N_i$ and $D_i$ preserve a uniform precision such that all the $n_i$ and $d_i$ share the same range $[0,ulp)$. According to (\ref{eqn:accumulative_error}), AAET after iteration 3 is $[-6ulp, 4ulp]$, which means that the precision loss after iteration 3 is beyond $2$ bits. However, it is possible to limit the total error to be less than $3$ bits by constraining the magnitude of CET. As mentioned above, CET will converge to a extremely small value after iteration 3 so that the above-mentioned requirement for it can be satisfied naturally. Therefore, to meet the precision requirement of $2^{-64}$ for the EP case, at least $3$ extra bits of precision for all the $N_i$ and $D_i$ is necessary if they preserve a uniform precision, i.e. the uniform $ulp$ is set to $2^{-67}$. For negative example, error will happen as sure as fate if the uniform $ulp$ becomes $2^{-66}$, which will also be verified along with the positive example.

Our discussion is based on the assumption of $ulp=2^{-67}$. According to (\ref{eqn:accumulative_error}) and (\ref{eqn:convergence_error}), the total approximate error after iteration 2 is
\begin{equation*}
    -Q\epsilon_1^{\prime 2} +\sum_{i=0}^{1}{Qd_i}-\sum_{i=0}^{2}{n_i}-Q\cdot f_{1},
\end{equation*}
which implies that the total error mainly consists of CET and the sub-term that contains $f_1$ of AET. Note that the sum of the precise analytical expressions for the aforementioned two terms is actually
\begin{equation*}
     -Q\epsilon_2=-Q\epsilon_1^{\prime 2}-Qf_1(1-\epsilon_1^{\prime}).
\end{equation*}
To meet the precision requirement of $2^{-53}$ for DP, the absolute value of the sum of $\epsilon_1^{\prime 2}$ and $f_1$ must be less than $2^{-54}$ if $Q$ approaches its upper bound $2$.

\begin{figure}
    \centering
    \includegraphics[width=\linewidth]{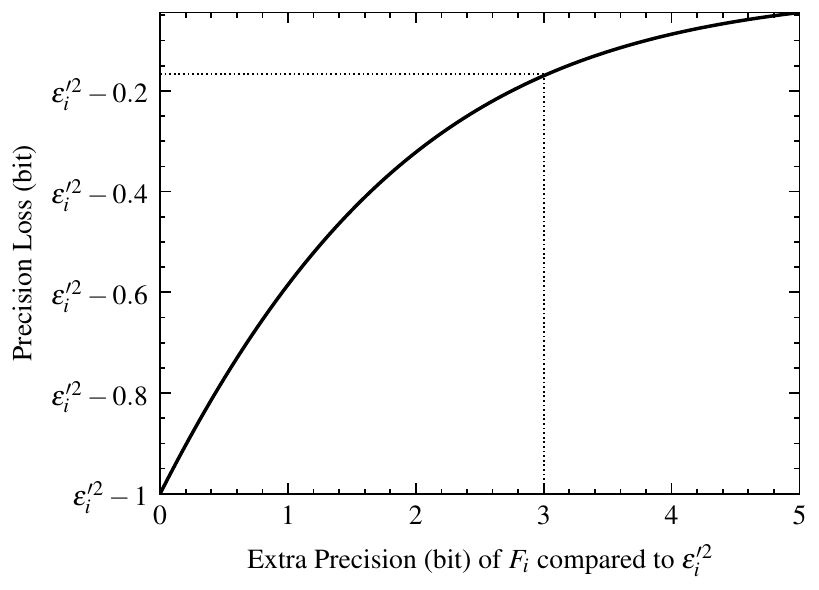}
    \caption{The precision loss corresponding to $F_i$ with different precision.}
    \label{fig:extra_precision}
\end{figure}
According to \figurename \ref{fig:extra_precision}, $3$ extra bits of precision for $F_1$ is appropriate with only about $0.2$ bits of precision loss. For this reason, the value of $\epsilon_1^{\prime 2}$ is assumed to approximately equal to $2^{-54.2}$ such that $\epsilon_1^{\prime}$ equals to $2^{-27.1}$. Similarly, the value of $\epsilon_0^{\prime}$ can be assumed to approximately equal to $2^{-\frac{(27.1+0.2)}{2}}=2^{-13.65}$. Moreover, since $\epsilon_0^{\prime}=\epsilon_0+d_0$, $d_0$ is too small to influence this equation significantly, the value $2^{-13.65}$ can be regarded as a upper bound of the relative error when selecting an appropriate look-up table configuration. Based on our experimental data, the bipartite table consisting of a large table of size $2^{9}\times 14$, a small table of size $2^{10}\times 6$ and a 14-bit subtractor with maximum absolute value of relative error for about $2^{-13.662378}$ are selected to generate the initial reciprocal. With this practical data of $|\epsilon_0|_{max}$, the above-mentioned derivation process should be reversed and carried out again.

For iteration 0, there is
\begin{equation*}
    \lceil |\epsilon_0^{\prime}|\rceil=\lceil|\epsilon_0|\rceil+\lceil d_0\rceil\approx2^{-13.662378}.
\end{equation*}
For iteration 1, since $\lceil|\epsilon_0^{\prime}|\rceil^2\approx2^{-27.234756}$, $3$ extra bits of precision based on $27$ bits for $F_0$ is appropriate:
\begin{equation*}
    \lceil\epsilon_1^{\prime}\rceil=\lceil|\epsilon_0^{\prime}|\rceil^2+(1+\lceil|\epsilon_0^{\prime}|\rceil)\lceil f_0\rceil+\lceil d_1\rceil \approx2^{-27.114905},
\end{equation*}
which implies that both CET and AET are not the bottleneck of precision for the SP case. Besides, the derived upper bound for $\epsilon_1^{\prime}$ meets the above-made assumption of $\epsilon_1^{\prime}\leq 2^{-27.1}$. 
For iteration 2, $\lceil\epsilon_1^{\prime}\rceil^2\approx2^{-54.229811}$, similarly, $3$ extra bits of precision based on $54$ bits for $F_1$ is appropriate:
\begin{equation*}
    \lceil\epsilon_2\rceil=\lceil\epsilon_1^{\prime}\rceil^2+(1+\lceil\epsilon_1^{\prime}\rceil)\lceil f_1\rceil\approx2^{-54.032467},
\end{equation*}
which satisfies the above-mentioned precision requirement. In accordance with the EP case, a negative example for the DP case will be given. If $F_1$ preserves only $2$ extra bits of precision, $\lceil\epsilon_2\rceil$ will increase to about $2^{-53.858898}$ such that error will definitely happen. Both of the positive and negative examples will be verified as well. Additionally, for the reason that both types of error terms are not the precision bottleneck in the SP case, only positive example will be tested.
\begin{figure}[t]
    \centering
    \includegraphics[width=\linewidth]{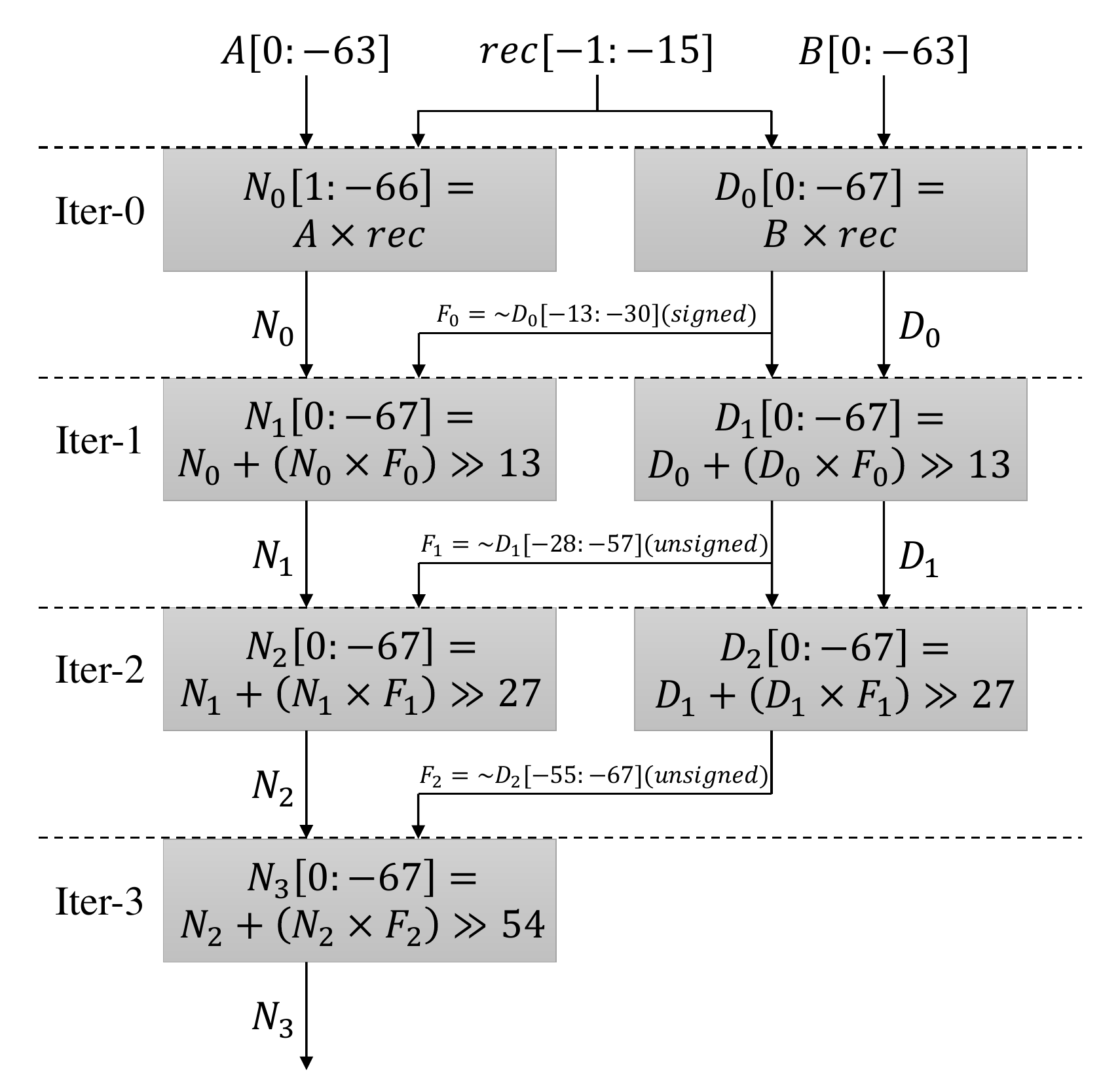}
    \caption{The iteration flow for the 3-stage divider, where $[a,b]$ represents the binary digits from $a$-th to $b$-th bit of the infinite precision format. Specifically, the first bits of the integer and fraction parts are numbered as the $0$-th and $-1$-th bit respectively.}
    \label{fig:3-staged}
\end{figure}


The above-mentioned analysis indicates the specific data-width configuration in this divider. However, it is worth noting that the sign of $\epsilon_0$ is not determined. If $\epsilon_0<0$ and $Q$ approaches to $2$, it is possible for $N_0$ to be greater than $2$. In this case, $1$ extra bit is necessary in the integer part of $N_0$, which leads to a mismatch in data-width of the multiplier if this divider is implemented by multiplier reusing. To solve this problem, design of $1$ less bit of precision only for $N_0$ can be considered, the doubled $n_0$ error in $N_0$ only causes the lower bound of AAET in the EP case to change from $-6ulp$ to $-7ulp$, the original design of $3$ extra bits for any other $N_i$ and $D_i$ still works.

\subsection{Two-Stage Divider for SP, DP and EP}\label{sec:stage-2}
In this subsection, a case that both the convergent and the accumulative error terms are the bottleneck of precision will be discussed. Consider a Goldschmidt divider with two iterations, where the result of the SP case is generated after iteration 1 and the results of DP and EP are available after iteration 2.

To start with, if the error requirement for the EP case can be met, then the requirement for the DP case will be satisfied affirmatively. The approximate accumulative error term after iteration 2 is $-5ulp$ to $+2ulp$, with a span of $7ulp$. According to \cite{kong2010rounding}, this asymmetrical error bound can be adjusted by adding a bias to decease its maximum absolute value, i.e. the final error can be constrained within $[-4ulp, 4ulp]$ such that $2$ extra bits of precision is enough. Nevertheless, considering that if $1$ bit less of precision is applied to $N_0$ as that in Section \ref{sec:stage-3}, the error redundancy remained for CET is almost $0$. Consequently, the extra precision is still $3$ bits here and a bias of $+5ulp$ is applied to the approximate quotient from iteration 2 to turn the range of AAET to $[-1ulp,+7ulp]$. As a result, there will be about $-7ulp$ of redundancy for CET.
\begin{figure}[t]
    \centering
    \includegraphics[width=\linewidth]{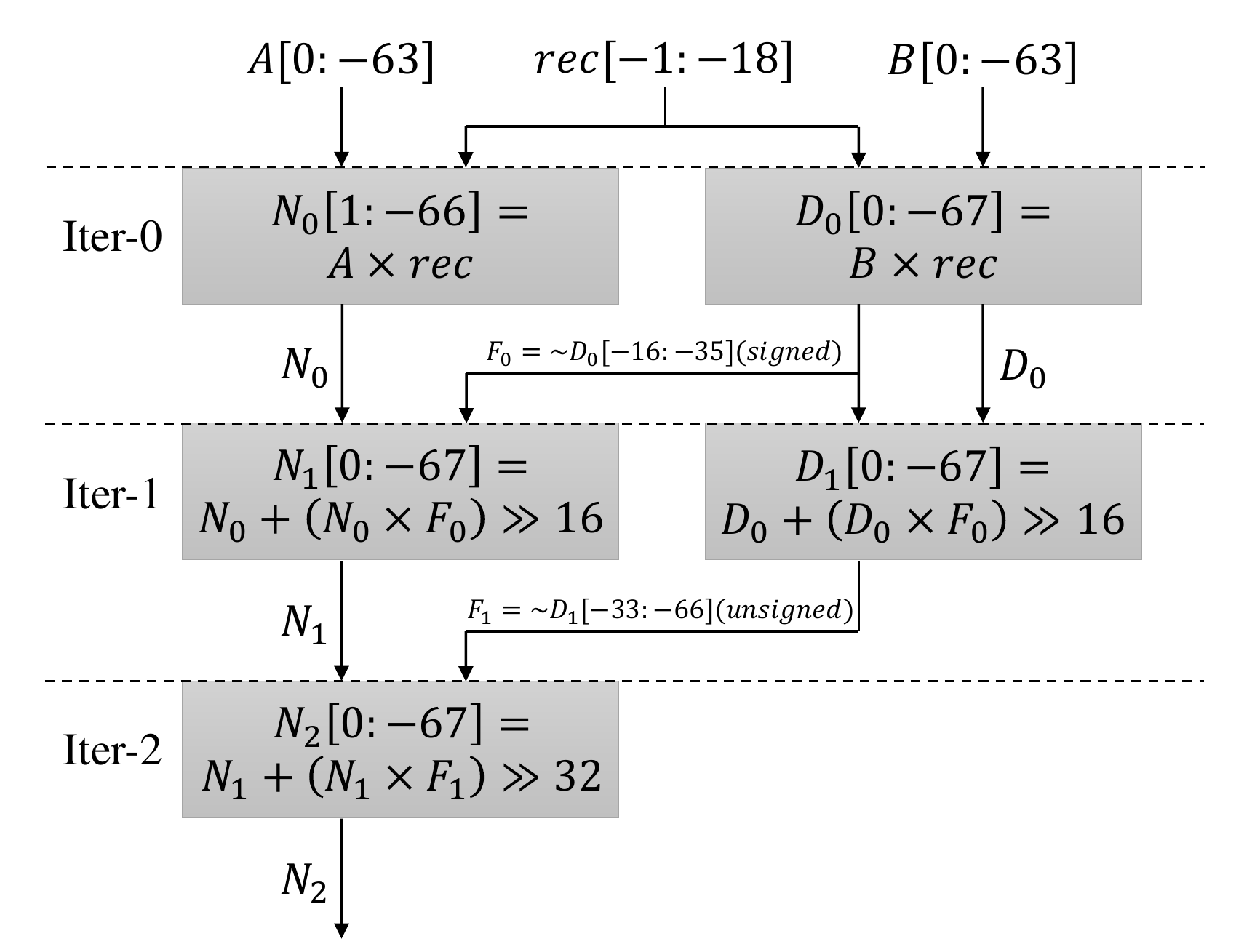}
    \caption{The iteration flow for the 2-stage divider.}
    \label{fig:2-staged}
\end{figure}

When $Q$ approaches $2$, $\epsilon_1^{\prime 2}$ should be less than $3.5\times2^{-67}$, i.e. $\epsilon_1^{\prime}<2^{-32.596322}$. The same as the train of thought in Section \ref{sec:stage-3}, $2^{-16.398161}$ can be regarded as a reference value to determine the target look-up table. In this example, the bipartite table consists of a large table of size $2^{11}\times 17$, a small table of size $2^{12}\times 7$ and a 17-bit subtractor with maximum absolute value of relative error for about $2^{-16.576687}$. For iteration 0, we have
\begin{equation*}
    \lceil |\epsilon_0^{\prime}|\rceil=\lceil|\epsilon_0|\rceil+\lceil d_0\rceil\approx2^{-16.576687}.
\end{equation*}
For iteration 1, since $\lceil|\epsilon_0^{\prime}|\rceil^2\approx2^{-33.1533754}$, only $2$ extra bits of precision based on $33$ bits for $F_0$ is enough:
\begin{equation*}
    \lceil\epsilon_1^{\prime}\rceil=\lceil|\epsilon_0^{\prime}|\rceil^2+(1+\lceil|\epsilon_0^{\prime}|\rceil)\lceil f_0\rceil+\lceil d_1\rceil \approx2^{-32.79972341}. 
\end{equation*}
This calculation result meets the precision requirements for $\epsilon_1^{\prime}$, and implies that neither of the two types of error are the bottleneck of precision for the SP case. Under this condition, the value of $\lfloor\text{CET}\rfloor\approx -5.280056\cdot ulp$ and the rough lower bound of the final error is then $-6.280056ulp$, which is far from $-8ulp$. Therefore, there is no need to do the careful numerical calculation mentioned in Section \ref{sec:accumulative_error}. While for a negative example, if there is only $1$ extra bit of precision for $F_0$, $\epsilon_1^{\prime}$ increases to about $2^{-32.515685}$ and $\lfloor\text{CET}\rfloor$ changes to $-7.827922ulp$, which is beyond the limit of $-7ulp$. Therefore, error will certainly appear. Noting that the precision of $F_{last}$ can be decreased by $1$ more bit, the rough lower bound in the positive example is only $-6.280056ulp$. If the value of $\lceil f_{last}\rceil$ doubles and the original rough lower bound turns to $-7.280056ulp$, which is still far from $-8ulp$. The new positive example and the negative one will be verified in the following section.

\section{Discussion}\label{sec:discussion}
In this section, the criticality problem in the proposed analysis theory which may lead to an uncertainty of data-width choice is analyzed. The verification results of the above-mentioned two-sided cases and the verifiability of random test are presented and discussed.

\begin{table}[t]
   \renewcommand{\arraystretch}{1.2}
  \setlength{\tabcolsep}{6pt}
    \centering
    \caption{The equivalent extra precision for the numerator and denominator for different methods.}
    \label{tab:equivalent_extra_precision}
    \begin{tabular}{|c|c|c|c|c|c|}
    \hline
       Method &  Proposed & \cite{kong2010rounding} & \cite{even2005parametric} & \cite{oberman1999floating}&\cite{schulte2007floating}\\
    \hline
    Extra precision & 3 bits&3 bits & 4 bits & 6 bits& 6 bits\\
    \hline
    \end{tabular}
\end{table}

\subsection{Criticality Problem}\label{sec:criticality}
To illustrate this problem, consider a case where the upper bound of its AAET is $+5ulp$ and the precision loss is obviously beyond $2$ bits so that the extra data-width has to be $3$ bits. However, if the upper bound of AAET is exactly $+4ulp$, $2$ or $3$ extra bits becomes a problem when pursuing for the optimal data-width.

According to the analysis in Section \ref{sec:accumulative_error}, the iterative factors $F_i$ with values slightly greater than $1$ will lead to a slightly greater upper bound of AET than $+4ulp$. Nevertheless, the $n_i$ and $d_i$ errors are rigorously less than $ulp$ at least, under the assumption of truncation in this paper. Assume that $F_i$ and another multiplicand have a fractional precision of $f$ and $k$ bits, respectively, then the upper bound is given by
\begin{equation*}
    2^{-k}-2^{-(k+f)}=(1-2^{-f})\times ulp.
\end{equation*}
As a result, it is necessary to analyze if the rigorous upper bound of AET can be less than $+4ulp$ to use only 2 extra bits.

Firstly, assume that $\epsilon_{i-1}=2^{-e}$ before an iteration begins. According to the data-width determination principle for $F_i$ in Section \ref{sec:stage-3} and \ref{sec:stage-2}, the length of the fractional part of $F_{i-1}$ will be about $\lfloor 2e+3 \rfloor$ bits. As a result, the upper bound for the corresponding $n_{i}$ and $d_{i}$ will be $(1-2^{-\lfloor 2e+3 \rfloor})\times ulp$ while the new $\epsilon_{i}$ will be greater than $2^{-2e}$ because of the precision loss. Then the coefficient $(1-2^{-\lfloor 2e+3 \rfloor})$ can be canceled by the new $F_i\approx 1+2^{-2e}$ in the next iteration.

Besides, there are still other elements may lead to the problem of "being less than $ulp$", e.g., the infinite precision quotient $Q$ which is rigorous less than $2$, the $n_{last}$ which is not multiplied by any iterative factors and the accumulative error cause by $f_{last}$. However, these coefficient of "being less than $ulp$" can be canceled easily by the multiplications with the first iterative factor $F_0$ with relatively greater $\epsilon_0$. Consequently, it is possible for the rigorous upper bound of the accumulative error term to be greater than $+4ulp$ and 3 extra bits of precision is indeed necessary.
\subsection{Verification Results}
\begin{figure}[t]
    \centering
    \includegraphics{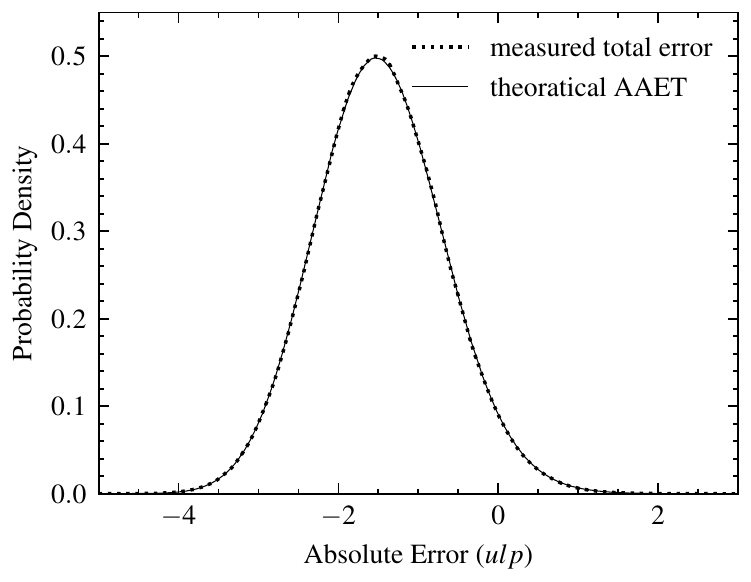}
    \caption{The fitted PDF of the total error measured in the 3rd stage of the divider in Section \ref{sec:stage-3} and its theoretical approximate PDF derived by AAET of this stage.}
    \label{fig:error_PDF}
\end{figure}
All of the positive examples mentioned in Section \ref{sec:stage-3} and \ref{sec:stage-2} have been verified by more than 10 billion of random test vectors. While for the negative examples, the divisors are set in a specific range to get the approximate reciprocal with the greatest initial iterative error to increase the efficiency of finding the inputs which lead to the incorrect rounding.

There is a problem of verifiability for the random test. The non-verifiable cases happen when the actual value of AET is required to approach the bounds of AAET. For example, if the added bias increases by $1ulp$ in Section \ref{sec:stage-2} and the upper bound of AAET after iteration 2 becomes $+8ulp$ from $+7ulp$. The analysis in Section \ref{sec:criticality} tells rounding error will happen but it is of great difficulty to proof that by recording the maximum positive error during the random test, since the probability is extremely small to encounter a case where all $d_i$ error approaches its maximum and all $n_i$ error approaches $0$. That is the reason why some of the negative examples are constructed by increasing the uniform $ulp$ instead of decreasing the precision of certain intermediate iterative results to make the bound of error lie on those critical values. In this view, the proposed analysis theory is able to determine whether a certain data-width is the optimal choice.
Moreover, for a practical example, the lower bound of AAET is $-6ulp$ in the case of Section \ref{sec:stage-3}, and according to the random test, the actual minimum negative error ever recorded after test vectors of more than 10 billion is about $-5.75ulp$.

\section{Conclusion}
In this paper, an error analysis theory for Goldschmidt division is proposed, which provides a clear view of how intermediate errors affect the final iterative result. It can not only work as a constraint in the optimization problem of data-width selection in Goldschmidt divider design, but also gives the basic thoughts of how to determine these data-widths. The proposed theory enables a state-of-the-art data-width reduction for the numerators and denominators, and leads to a significant data-width reduction for the iterative factors through numerical method with detailed theoretical basis. It can also be used in Goldschmidt division to enable a faster than quadratic convergence. This error analysis has been verified by 100 billion of random test vectors on two divider models implemented by VerilogHDL.

\ifCLASSOPTIONcompsoc



\ifCLASSOPTIONcaptionsoff
  \newpage
\fi



\bibliographystyle{IEEEtran}
\bibliography{IEEEabrv,ref.bib}

\begin{thebibliography}{10}
\providecommand{\url}[1]{#1}
\csname url@samestyle\endcsname
\providecommand{\newblock}{\relax}
\providecommand{\bibinfo}[2]{#2}
\providecommand{\BIBentrySTDinterwordspacing}{\spaceskip=0pt\relax}
\providecommand{\BIBentryALTinterwordstretchfactor}{4}
\providecommand{\BIBentryALTinterwordspacing}{\spaceskip=\fontdimen2\font plus
\BIBentryALTinterwordstretchfactor\fontdimen3\font minus
  \fontdimen4\font\relax}
\providecommand{\BIBforeignlanguage}[2]{{%
\expandafter\ifx\csname l@#1\endcsname\relax
\typeout{** WARNING: IEEEtran.bst: No hyphenation pattern has been}%
\typeout{** loaded for the language `#1'. Using the pattern for}%
\typeout{** the default language instead.}%
\else
\language=\csname l@#1\endcsname
\fi
#2}}
\providecommand{\BIBdecl}{\relax}
\BIBdecl

\bibitem{goldschmidt1964applications}
R.~E. Goldschmidt, ``Applications of division by convergence,'' Ph.D.
  dissertation, Massachusetts Institute of Technology, 1964.

\bibitem{parhami2010computer}
B.~Parhami, \emph{Computer arithmetic}.\hskip 1em plus 0.5em minus 0.4em\relax
  Oxford university press, 2010, vol.~20, no.~00.

\bibitem{beame1986log}
P.~W. Beame, S.~A. Cook, and H.~J. Hoover, ``Log depth circuits for division
  and related problems,'' \emph{SIAM Journal on Computing}, vol.~15, no.~4, pp.
  994--1003, 1986.

\bibitem{krishnamurthy1970optimal}
E.~Krishnamurthy, ``On optimal ierative schemes for high-speed division,''
  \emph{IEEE Transactions on Computers}, vol. 100, no.~3, pp. 227--231, 1970.

\bibitem{oberman1999floating}
S.~F. Oberman, ``Floating point division and square root algorithms and
  implementation in the {AMD-K7}/sup {TM}/microprocessor,'' in
  \emph{Proceedings 14th IEEE Symposium on Computer Arithmetic}.\hskip 1em plus
  0.5em minus 0.4em\relax IEEE, 1999, pp. 106--115.

\bibitem{schwarz1997cmos}
E.~M. Schwarz, L.~Sigal, and T.~J. McPherson, ``{CMOS} floating-point unit for
  the {S/390} parallel enterprise server {G4},'' \emph{IBM Journal of Research
  and Development}, vol.~41, no. 4.5, pp. 475--488, 1997.

\bibitem{russinoff1998mechanically}
D.~M. Russinoff, ``A mechanically checked proof of ieee compliance of the
  floating point multiplication, division and square root algorithms of the
  amd-k7™ processor,'' \emph{LMS Journal of Computation and Mathematics},
  vol.~1, pp. 148--200, 1998.

\bibitem{even2005parametric}
G.~Even, P.-M. Seidel, and W.~E. Ferguson, ``A parametric error analysis of
  goldschmidt's division algorithm,'' \emph{Journal of Computer and System
  Sciences}, vol.~70, no.~1, pp. 118--139, 2005.

\bibitem{kong2010rounding}
I.~Kong and E.~E. Swartzlander, ``A rounding method to reduce the required
  multiplier precision for {Goldschmidt} division,'' \emph{IEEE Transactions on
  Computers}, vol.~59, no.~12, pp. 1703--1708, 2010.

\bibitem{schulte2007floating}
M.~J. Schulte, D.~Tan, and C.~E. Lemonds, ``Floating-point division algorithms
  for an x86 microprocessor with a rectangular multiplier,'' in \emph{2007 25th
  International Conference on Computer Design}.\hskip 1em plus 0.5em minus
  0.4em\relax IEEE, 2007, pp. 304--310.

\bibitem{kong2009goldschmidt}
I.~Kong and E.~E. Swartzlander, ``A goldschmidt division method with faster
  than quadratic convergence,'' \emph{IEEE Transactions on Very Large Scale
  Integration Systems}, vol.~19, no.~4, pp. 696--700, 2009.

\bibitem{sarma1995faithful}
D.~D. Sarma and D.~W. Matula, ``Faithful bipartite rom reciprocal tables,'' in
  \emph{Proceedings of the 12th Symposium on Computer Arithmetic}.\hskip 1em
  plus 0.5em minus 0.4em\relax IEEE, 1995, pp. 17--28.

\bibitem{schulte1997symmetric}
M.~J. Schulte and J.~E. Stine, ``Symmetric bipartite tables for accurate
  function approximation,'' in \emph{Proceedings 13th IEEE Sympsoium on
  Computer Arithmetic}.\hskip 1em plus 0.5em minus 0.4em\relax IEEE, 1997, pp.
  175--183.

\bibitem{farmwald1981design}
P.~M. Farmwald, ``On the design of high performance digital arithmetic units,''
  Ph.D. dissertation, Stanford University, 1981.

\bibitem{nakano1987method}
H.~Nakano, ``Method and apparatus for division using interpolation
  approximation,'' Nov.~17 1987, {US} Patent 4, 707, 798.

\bibitem{piso2010variable}
D.~Piso and J.~D. Bruguera, ``Variable latency goldschmidt algorithm based on a
  new rounding method and a remainder estimate,'' \emph{IEEE Transactions on
  Computers}, vol.~60, no.~11, pp. 1535--1546, 2010.

\bibitem{obermann1997division}
S.~F. Obermann and M.~J. Flynn, ``Division algorithms and implementations,''
  \emph{IEEE Transactions on Computers}, vol.~46, no.~8, pp. 833--854, 1997.

\bibitem{cornea1999correctness}
M.~A. Cornea-Hasegan, R.~A. Golliver, and P.~Markstein, ``Correctness proofs
  outline for newton-raphson based floating-point divide and square root
  algorithms,'' in \emph{Proceedings 14th IEEE Symposium on Computer
  Arithmetic}.\hskip 1em plus 0.5em minus 0.4em\relax IEEE, 1999, pp. 96--105.

\end{thebibliography}
%



%

\begin{IEEEbiography}[{\includegraphics[width=1.05in,height=1.31in,clip,keepaspectratio]{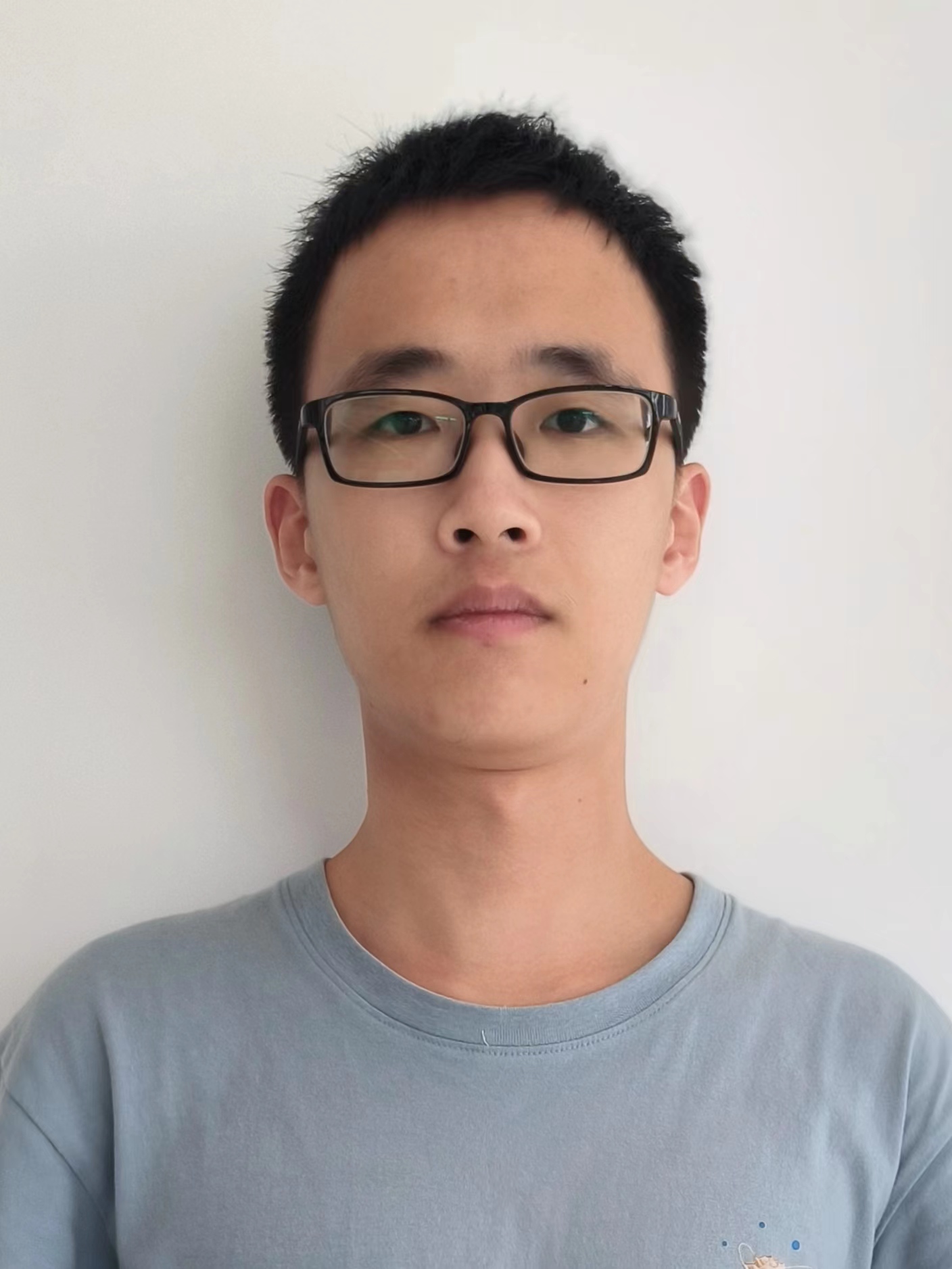}}]{Binzhe Yuan} received the B.Eng. degree in Electronic Information Engineering from ShanghaiTech University, Shanghai, China, in 2022. He is currently pursuing the M.Eng. degree in ShanghaiTech University, Shanghai, China. His current research interests include VLSI design and computer arithmetic.
\end{IEEEbiography}
\begin{IEEEbiography}[{\includegraphics[width=1.05in,height=1.31in,clip,keepaspectratio]{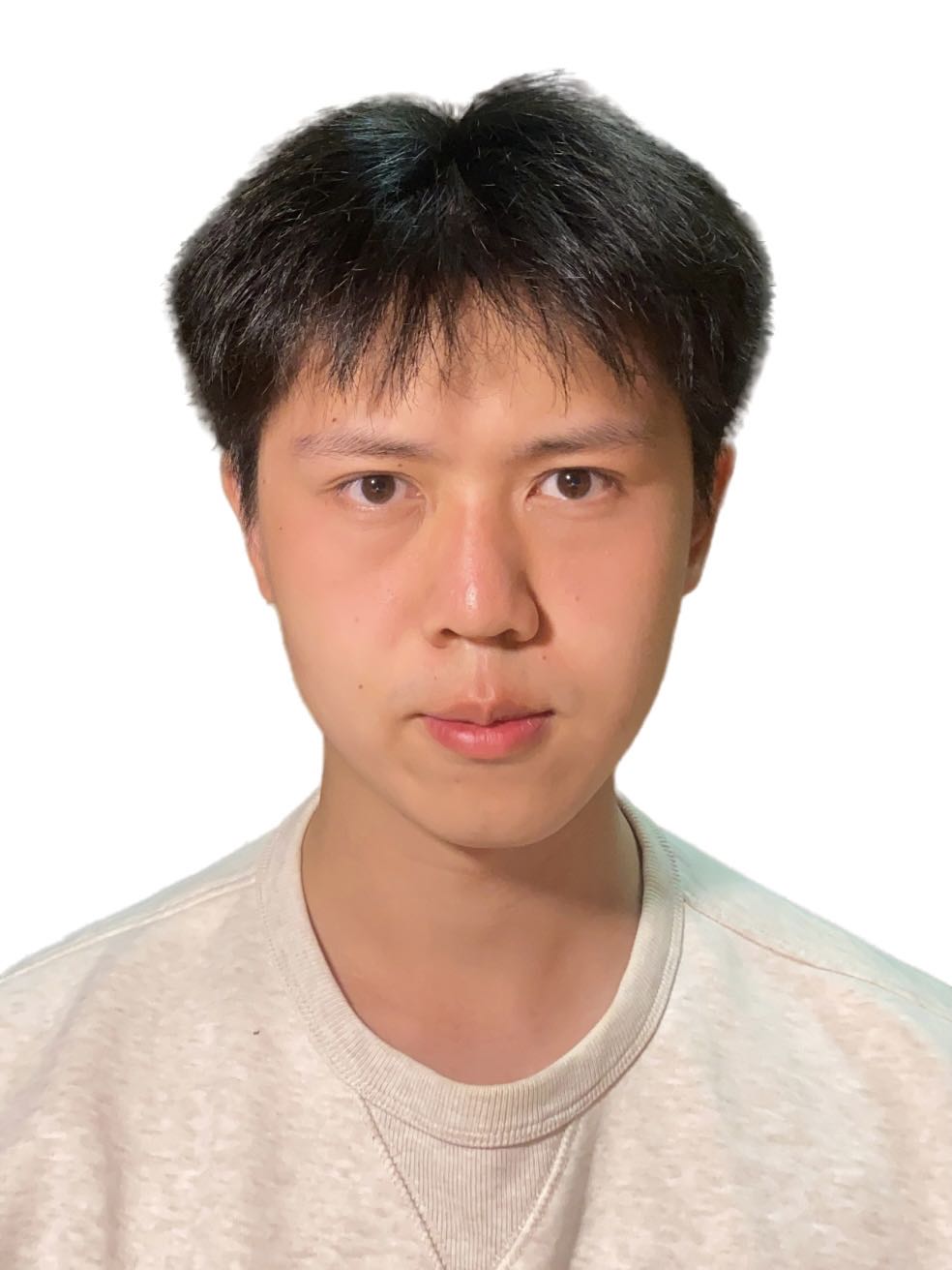}}]{Liangtao Dai} is currently pursuing the B.Eng. degree in ShanghaiTech University, Shanghai, China. His current research interest is VLSI design.
\end{IEEEbiography}
\begin{IEEEbiography}[{\includegraphics[width=1in,height=1.25in,clip,keepaspectratio]{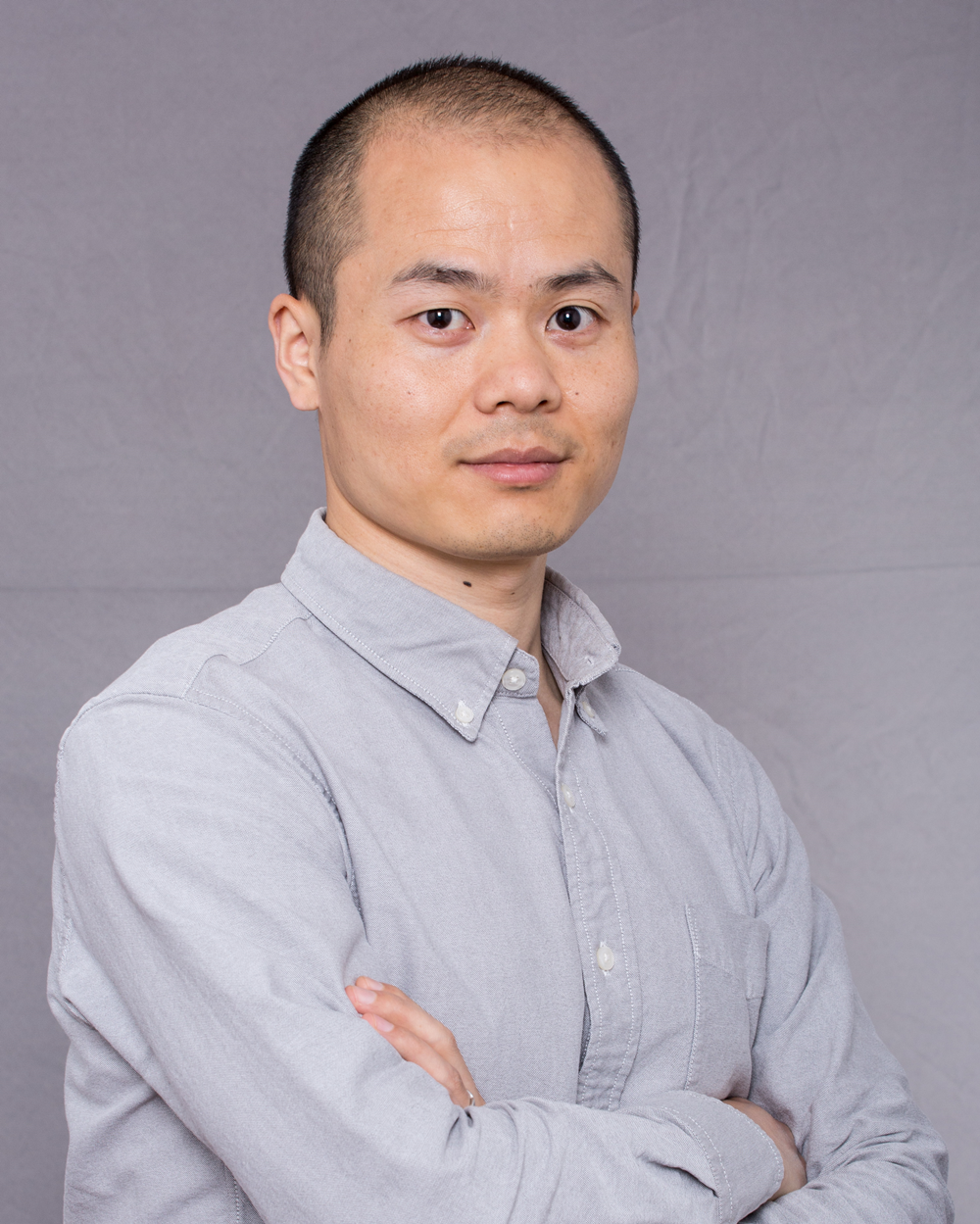}}]{Xin Lou}(Senior Member, IEEE) received the B.Eng. degree in Electronic Information Technology and Instrumentation from Zhejiang University (ZJU), China, in 2010 and M.Sc. degree in System-on-Chip Design from Royal Institute of Technology (KTH), Sweden, in 2012 and PhD degree in Electrical and Electronic Engineering from Nanyang Technological University (NTU), Singapore, in 2016. Then he joined VIRTUS, IC Design Centre of Excellence at NTU as a research scientist. In Mar. 2017, he joined the School of Information Science and Technology, ShanghaiTech University, Shanghai, China, as an assistant professor and establish the VLSI Signal Processing Lab. His research interests include energy-efficient integrated circuits for machine vision and 3D graphics, and other VLSI digital signal processing systems.
\end{IEEEbiography}
\vfill




\end{document}